\documentclass[12pt]{amsart}
\usepackage{amsfonts,amssymb,amsmath,textcomp}

 2

\newtheorem{thm}{Theorem}[section]

\newcommand{\thmref}[1]{Theorem~\ref{#1}}

\begin{document}
\title
[Sign changes]{A short note on Sign changes}

\author[J. Meher]{Jaban Meher}
\address[Jaban Meher]{Department of Mathematics, Queen's University\\
Kingston, Ontario, K7L 3N6, Canada}
\email[Jaban Meher]{jaban@mast.queensu.ca}

\author[K. D. Shankhadhar]{Karam Deo Shankhadhar }
\address[Karam Deo Shankhadhar]{Harish-Chandra Research Institute\\
Chhatnag Road, Jhunsi\\
Allahabad 211019, India}
\email[Karam Deo Shankhadhar]{karam@hri.res.in}

\author{G.K. Viswanadham }
\address[G.K. Viswanadham ]{Harish-Chandra Research Institute\\
Chhatnag Road, Jhunsi\\
Allahabad 211019, India}
\email[G.K. Viswanadham]{viswanath@hri.res.in}

\maketitle 

\vspace{.2cm}

\begin{abstract}
In this paper, we present a quantitative result for the number of sign changes for the sequences
$\{a(n^j)\}_{n\ge 1}, j=2,3,4$ of the Fourier coefficients of normalized Hecke eigen cusp forms 
for the full modular group $SL_2(\mathbb{Z})$. We also prove a similar kind of quantitative result
for the number of sign changes of the $q$-exponents 
$c(p)~(p \mbox{~vary ~over ~primes})$ of certain generalized modular functions for the 
congruence subgroup $\Gamma_0(N)$, where $N$ is square-free. 
\end {abstract}

\vspace{.1cm}
\section{Introduction}
\vspace{.1cm}
The sign changes of the Fourier coefficients of modular forms are well studied. It is known that,
if the Fourier coefficients of a cusp form are real, then they change signs infinitely 
often. The proof of this fact follows from the Landau's theorem on Dirichlet series with 
non-negative coefficients. Further, many quantitative results for the number of sign changes for 
the sequence of the Fourier coefficients have been established. The sign changes of the subsequence
of the Fourier coefficients at prime numbers
was  first studied by M. Ram Murty \cite{murty}. 

In our first result, 
we show that if $h(z)=\sum_{n\ge 1}a(n)q^n$ is a Hecke eigenform
for the full modular group $SL_2(\mathbb{Z})$, then
the subsequences $\{a(n^j)\}_{n\ge 1}, j=2, 3, 4$ of the Fourier coefficients change signs infinitely often. 
In fact, we give a lower 
bound for the number of sign changes in the interval $(x,2x]$. We use the results of O. M. Fomenko \cite{F},
H. Lao and A. Sankaranarayanan \cite{LS} and G. S. L\"u \cite{LU} to prove this result.

Our second result is about the sign changes of the $q$-exponents of a certain class of generalized
modular functions $f$ on the Hecke congruence subgroup $\Gamma_0(N)$. We always assume that $f$ is non-constant 
and the weight of $f$ is zero. A generalized modular
function $f$ is a holomorphic function on the upper half plane $\mathcal{H}$, meromorphic at the
cusps, that transforms under $\Gamma_0(N)$ like a usual modular function, with the exception
that the character $\chi$ need not be unitary, but $\chi(\gamma)=1$ for all parabolic
$\gamma\in \Gamma_0(N)$ of trace $2$. For more details on generalized modular functions,
see \cite{KNM} and \cite{KOM}.
It is known that, every generalized modular function has a product expansion
$$
f(z)=cq^h\prod{(1-q^n)}^{c(n)} ~~~~~~(0<|q|<\epsilon),
$$
where $h\in \mathbb{Z}$, $c$ and $c(n)$ are uniquely determined complex numbers and
$q=e^{2\pi iz} ~(z\in \mathcal{H})$.
In \cite{KM}, it has been proved that if $f$ is a generalized modular function with
$div(f)=\emptyset$, then $c(n)~(n\ge 1)$ change signs infinitely often. This result was
proved by using a theorem of \cite{KNM} that any generalized modular function
$f$ with $div(f)=\emptyset$ corresponds to a cusp form of weight $2$ by taking the
logarithmic derivative. However, the result is not quantitative. Using the results of 
Y. J. Choie and W. Kohnen \cite{choie-kohnen}, H. Iwaniec, W. Kohnen and J. Sengupta \cite{i-k-s} 
and the 
relation between the $q$-exponents of generalized modular function and the Fourier
coefficients of the corresponding cusp form of weight $2$, one can deduce a bound for the first sign 
change of the $q$-exponents of a certain kind of generalized modular function.

In our second result, 
we give a quantitative result for the number of sign changes of the subsequence of $q$-exponents at primes,
$\{c(p)\}$, in the interval $(x, 2x]$ for 
a certain kind of generalized modular function on $\Gamma_0(N)$, where $N$ is square-free.
To prove this result, we use the method given by Y. J. Choie, A. Sankaranarayanan and J. Sengupta
\cite{css}. 
 
\section{Statement of results}\label{results}

Let $h(z)=\sum_{n\ge 1}a(n)q^n$ be a Hecke eigenform of weight $k$ on the
full modular group $SL_2(\mathbb{Z})$. Let $\lambda(n)=\frac{a(n)}{n^{(k-1)/2}}$.
Denote by $\delta_j=2/11, 1/9, 2/27$ for $j=2, 3, 4$ respectively. Now, we state
our first result.
\begin{thm}\label{thm2}
For any $j\in\{2,3,4\}$, the sequence $\{\lambda(n^j)\}_{n\ge 1}$ change signs infinitely often. Moreover, the
sequence has at least $>>x^{\delta_j-2\epsilon}$ sign changes in the interval $(x,2x]$ 
for sufficiently large $x$, where $\epsilon$ is any small positive constant.
\end{thm}
From \cite[Theorem 2]{KNM}, it is known that if $f$ is a generalized modular function on $\Gamma_0(N)$ with $div(f)=\emptyset$,
then its logarithmic derivative $g=\frac{1}{2\pi i}\frac{f'}{f}$ is a cusp form of weight $2$
on $\Gamma_0(N)$ with trivial character. We now state our second result.
\begin{thm}\label{thm1}
Let $N$ be a square-free positive integer and $f$ be a non-constant generalized
modular function on $\Gamma_0(N)$ with $div(f)=\emptyset$. Suppose that the 
logarithmic derivative $g=\frac{1}{2\pi i}\frac{f'}{f}$ is a normalized newform. Then the 
sequence $c(p) (p ~\rm{prime})$ change signs infinitely often and it has 
at least $>>e^{A\sqrt{\log{x}}}$ sign changes in the interval $(x,2x]$ for sufficiently large 
$x$, where $A$ is an absolute constant.
\end{thm}

\section{Proofs} 

\noindent{\bf Proof of \thmref{thm2}}
Let us denote $1/2,3/4$ and $7/9$ by $\beta_j$ for $j=2,3$ and $4$ 
respectively. For any $\epsilon>0$, \cite{F} and \cite{LU} give the following 
estimates. 
\begin{equation}\label{eq:10}
\sum_{n\leq x}\lambda(n^j)\ll_{f,\epsilon} x^{\beta_j+\epsilon}.
\end{equation}
From \cite[Theorem 1.1, 1.2, 1.3]{LS}, we get the following estimates.
\begin{equation}\label{eq:11}
\sum_{n\leq x}\lambda^2(n^j)=B_jx+O_{f,\epsilon}(x^{1-\delta_j+\epsilon}),
\end{equation}
where $B_j$ are absolute constants, $\delta_j$ are defined as before and
these estimates are valid for any $\epsilon >0$. Let 
$h=h(x)=x^{1-\delta_j+2\epsilon}$, where $\epsilon$ is sufficiently small. 
We assume that the sequence 
$\{\lambda(n^j)\}_{n\geq 1}$ are of constant sign say positive for all 
$n\in (x,x+h]$. Using \eqref{eq:11}, we get 
\begin{equation}\label{eq:12}
\sum_{x<n\leq x+h}\lambda^2(n^j)=B_jh+O_{f,\epsilon}(x^{1-\delta_j+\epsilon})
\gg x^{1-\delta_j+2\epsilon}.
\end{equation}
On the other hand, using \eqref{eq:10}, we get 
\begin{equation}\label{eq:13}
\begin{split}
\sum_{x<n\leq x+h}\lambda^2(n^j)=\displaystyle\sum_{x<n\leq x+h}\lambda(n^j)
\lambda(n^j) \ll x^{2\epsilon}\displaystyle\sum_{x<n\leq x+h}
\lambda(n^j) & \ll x^{2\epsilon}\left((x+h)^{\beta_j+\epsilon}
+x^{\beta_j+\epsilon}\right)\\
& \ll x^{\beta_j+3\epsilon}.
\end{split}
\end{equation}
Now comparing $1-\delta_j$ and $\beta_j$ for $j=2,3,4$, we see that the bounds 
in \eqref{eq:12} and 
\eqref{eq:13} for $\displaystyle\sum_{x<n\leq x+h}\lambda^2(n^j)$ contradict
each other. Therefore, at least one $\lambda(n^j)$ for $x<n\leq x+h$ 
must be negative. Hence the sparse sequences $\{\lambda(n^j)\}_{n\geq 1}$
for $j=2,3,4$ change signs infinitely often and there are at least 
$\gg x^{\delta_j-2\epsilon}$ sign changes in the interval $(x,2x]$.

\smallskip

\noindent{\bf Proof of \thmref{thm1}}
Since $div(f)=\emptyset$, we have $f(z)=\prod_{n\ge 1}{(1-q^n)}^{c(n)}$.
Let $g(z)=\sum_{n\ge 1}b(n)q^n$  and $\lambda(n)=\frac{b(n)}{n^{1/2}}$. Then we have
$$
b(n)=-\displaystyle\sum_ {d|n} dc(d).
$$      
In particular, $b(1)=1=-c(1)$ and $b(p)=1-pc(p)$  for any prime $p$. 
 
\noindent Using a theorem of C. J. Moreno (\cite{M}) and a result of 
J. Hoffstein and D. Ramakrishnan (\cite{HR}) about the nonexistence of Siegel
zero, we have the following estimate. 
\begin{equation}\label{eq:2}
\sum_{p\leq x}\lambda(p)\log p=O(xe^{-A_1\sqrt{\log x}}),
\end{equation} 
where $A_1$ is an absolute constant. Following the proof of \cite[Theorem 4.1]{S} and using a nonvanishing result for symmetric square $L$-function
of a newform on square-free level
due to D. Goldfeld,
J. Hoffstein and D. Liemann (\cite{GHL}) for the proof of \cite[Lemma 4.1]{S}, we get the following estimate.
\begin{equation}\label{eq:3}
\sum_{p\leq x}\lambda^2(p)\log p=x+O(xe^{-A_2\sqrt{\log x}}),
\end{equation} 
where $ A_2 $ is an absolute constant. 

\noindent Since $\lambda(p)=\frac{b(p)}
{\sqrt{p}}=-\left(\frac{-1}{\sqrt{p}}+\sqrt{p}c(p)\right)$, therefore 
$\sqrt{p}c(p)=
\frac{1}{\sqrt{p}}-\lambda(p)$.
Let $c'(p)=\sqrt{p}c(p)$ for all primes $p$. The behaviour of $c'(p)$ and 
$c(p)$ are same in the sense of sign changes since ~$\frac{c'(p)}{c(p)}$~
is a positive 
real number. We have
\begin{eqnarray*}
\displaystyle\sum_{p\leq x}c'(p)\log p & = & \displaystyle\sum_{p\leq x}
\left(\frac{1}{\sqrt{p}}-\lambda(p)\right)\log p\\
                                         & = & O(\sqrt{x}\log x ) + O(xe^{-A_1\sqrt{\log x}}),
\end{eqnarray*}
here we used \eqref{eq:2} for getting the last line. Hence
\begin{equation}\label{eq:4}
\displaystyle\sum_{p\leq x}c'(p)\log p= O(xe^{-A_3\sqrt{\log x}}),~~~~~~A_3\rm{~is 
~an~absolute~constant}.
\end{equation}
We have
\begin{eqnarray*}
\displaystyle\sum_{p\leq x}c'^2(p)\log p & = & \displaystyle\sum_{p\leq x}
\left(\frac{1}{\sqrt{p}}-\lambda(p)\right)^2\log p\\
& = & \displaystyle 
\sum_{p\leq x} \lambda^2(p)\log p +2\displaystyle\sum_{p\leq x}\frac{c'(p)\log p}
{p}-
\displaystyle \sum_{p\leq x} \frac{\log p}{p}.
\end{eqnarray*}
Now, estimating first and second terms by using \eqref{eq:3} and \eqref{eq:4} 
respectively, 
we get 
the following.
\begin{equation}\label{eq:5}
\displaystyle\sum_{p\leq x}c'^2(p)\log p= x+O(x e^{-A_4\sqrt{\log x}}),
~~~~~~A_4\rm{~is 
~an~absolute~constant}.
\end{equation}
Let $h=h(x)$ be any function of $x$ with the property that $0<h(x)<x$. 
We evaluate $\displaystyle\sum_{x<p\leq x+h}c'(p)$ 
by two different ways and arrive at a  contradiction if $c'(p)$
do not change the sign as stated in the result.
Assume that the $c'(p)$ are of constant sign for $x<p\leq x+h$. Without loss of
generality, we assume that $c'(p)\geq 0$ for $x<p\leq x+h$. Using the fact that
$|c'(p)| \leq 3$, we get the following.
\begin{equation}\label{eq:6}
\sum_{x<p\leq x+h}c'(p) \geq A_5 \sum_{x<p\leq x+h}c'^2(p),
\end{equation}
where $A_5$ is an absolute constant. Using the estimate \eqref{eq:5}, we have
\begin{eqnarray*}
\displaystyle\sum_{p\leq x}c'^2(p) & = & \displaystyle\sum_{p\leq x}c'^2(p)
\frac{\log p}{\log p}\\
                                    & = & \frac{1}{\log x}\left
(\displaystyle\sum_{p\leq x}c'^2(p)\log p\right)+\int_2^x\frac{\left(
\displaystyle\sum_{p\leq t}c'^2(p)\log p\right)}{t \log^2 t}dt\\
                                    & = & \frac{1}{\log x}
(x+O(xe^{-A_4\sqrt{\log x}})+\int_2^x\frac{t+O(te^{-A_4\sqrt{\log t }})}
{t\log^2 t}dt\\
                                    & = & \frac{x}{\log x}+O (\frac{x}{\log x}
e^{-A_4\sqrt{\log x}})+\int_2^x\frac{1}{\log^2t }dt +O\left(\int_2^x\frac{e^{-A_4
\sqrt{\log t}}}{\log^2t}dt\right).\\
\end{eqnarray*}
From the above, we get
\begin{equation}\label{eq:7}
\displaystyle\sum_{p\leq x}c'^2(p)\sim \frac{x}{\log x}.
\end{equation}
Combining \eqref{eq:6} and \eqref{eq:7}, we deduce
\begin{equation}\label{eq:8}
\sum_{x<p\leq x+h}c'(p)\geq A_5\left(\displaystyle\sum_{p\leq x+h}c'^2(p)-
\displaystyle\sum_{p\leq x}c'^2(p)\right) \geq A_6 \left(\frac{x+h}{\log (x+h)}
-\frac{x}{\log x}\right)
\gg \frac{h}{\log x},
\end{equation}
where $A_6$ is an absolute constant. 

\noindent On the other hand using \eqref{eq:4}, 
we have the following estimate.
\begin{eqnarray*}
\sum_{x<p\leq x+h}c'(p) & = & \displaystyle \sum_{x<p\leq x+h}c'(p)
\frac{\log p}{\log p}\\
                                      & \leq & \frac{1}{\log x}\displaystyle
\sum_{x<p\leq x+h}c'(p)\log p\\
                                      & \ll & \frac{1}{\log x} \left((x+h)
e^{-A_3\sqrt{\log(x+h)}}+xe^{-A_3\sqrt{\log x}}\right).
\end{eqnarray*}
From the above estimate, we obtain
\begin{equation}\label{eq:9}
\sum_{x<p\leq x+h}c'(p) \ll  \frac{x}{\log x}e^{-A_3\sqrt{\log x}}.
\end{equation}
Choosing an appropriate absolute constant $A$ and $h(x)=
\frac{x}{e^{A\sqrt{\log x}}}$, we see that the bounds obtained in 
\eqref{eq:8} and \eqref{eq:9} contradict each other. 
Hence, at least one $c'(p)$ for $x<p\leq x+h$ must be negative. This implies that 
the sequence $\{c(p)\}$ has infinitely many sign changes and there are at least 
$\gg e^{A\sqrt{\log x}}$ sign changes for the sequence $\{c(p)\}$, 
whenever $p\in(x,2x]$.
\\
\qed

\noindent{\bf Acknowledgments}\\
We thank Prof. A. Sankaranarayanan for useful discussion and giving valuable 
suggestions. Also, we thank the referee for correcting some errors.

\end{document}